%

\magnification=\magstep1
\input amstex
\documentstyle{amsppt}

\define\prg{\Cal G_{c\& c}}
\topmatter
\title
More on the cut and choose game 
\endtitle
\author
Jind\v rich Zapletal
\endauthor
\thanks
 \endthanks
\affil
The Pennsylvania State University
\endaffil
\address
Department of Mathematics, 
The Pennsylvania State University, 
University Park,  PA 16802
\endaddress
\email
zapletal\@math. psu. edu
\endemail
\subjclass
\endsubjclass
\abstract
We improve on the results of \cite {Vel} and give some examples.
\endabstract
\endtopmatter
\document
In \cite {Jech}, a number of infinite games on (complete) Boolean algebras 
were defined. Among them was the following Prikry-style ``cut and choose" 
game $\prg :$

\medpagebreak
\settabs \+\indent&aaaaa&bbbbb&ccc&ddd&eee&ratattatat&ggg&hhh&ratatttatat&brum&&\cr
\+&I&$p,b_0$&&$b_1$&&\dotfill&$b_i$&&\dotfill&&$i<\omega$&\cr
\+&II&&$r_0$&&$r_1$&\dotfill&&$r_i$&\dotfill&\cr
\medpagebreak

where $p,$ $b_i$'s are elements of the Boolean algebra $B$ in question 
and $r_i\in 2$ for $i<\omega .$ II wins a run of the game $\prg$ if 
$p\land \underset {i<\omega ,r_i=1}\to \bigwedge b_i -\underset 
{i<\omega ,r_i=0}\to \bigvee b_i\neq 0.$ Intuitively, I chooses an element 
$p\in B$ under which the game is to take place and then asks questions like 
``is $b_i$ an element of the generic ultrafilter?" II must give an answer. 
II wins if in the end the choices he has made are compatible with $p.$ 
For this game to be interesting it is necessary that $B$ be complete and 
we shall assume so for all Boolean algebras considered in the rest of the 
paper. \cite {Jech} gives the following:
\proclaim
{Fact 1} I has a winning strategy in $\prg$ if and only if $B$ adds 
real numbers.
\endproclaim
\proclaim
{Problem 1} Is it true in ZFC that II has a winning strategy in 
$\prg$ does not imply that $B$ has a $\sigma$-closed dense subset?
\endproclaim
\proclaim
{Problem 2} Find in ZFC a Boolean algebra $B$ such that $\prg$ is 
undetermined.
\endproclaim
In this paper, we give a solution to the Problem 2 and an approximation to 
a solution of Problem 1, best possible modulo the solution of 
the apparently difficult open Problem 6 of \cite {Jech}. (This has
been done in \cite {Vel} by different methods.) Further we 
consider relationship between $\prg$ and some stronger distributivity 
properties of $B.$ 

Let us fix a  Boolean algebra $B.$ A strategy for II in $\prg$ is 
a function $\sigma: B\times B^{<\omega }\to 2^{<\omega }$ such that 
if $\langle p,t\rangle \in dom(\sigma ),$ $lth(t)=i$ implies 
$lth(\sigma (\langle p,t\rangle ))=i$ and for each $j<i,$ 
$\sigma (\langle p,t\restriction j\rangle )=\sigma 
(\langle p,t\rangle )\restriction j.$ A strategy $\sigma$ is winning 
if for all $\langle p,t\rangle \in B\times B^\omega ,$
 $p\land \underset {\sigma (\langle p,t\restriction i+1\rangle )(i)=1}
\to \bigwedge t(i)-\underset {\sigma (\langle p,t\restriction i+1
\rangle )(i)=0}\to \bigvee t(i)\neq 0.$

There are two canonical examples of Boolean algebras in which II has 
a winning strategy in $\prg :$ namely, any $\sigma$-closed forcing and 
Prikry forcing with a normal measure $U$ \cite {P\v r\' \i}.
 (See Example 3 for a completely different poset.) 
We can paraphrase the Problem 1 as 
``Find an analog of Prikry forcing in ZFC." Our results are as follows:

\proclaim
{Theorem 1} If II has a winning strategy in $\prg ,$ then $B$ is semiproper.
\endproclaim

\proclaim
{Theorem 2}  If II has a winning strategy in $\prg ,$ and $B$ has 
$(2^\omega )^+$-c.c. then II has a winning strategy in the descending 
chain game $\Cal G.$
\endproclaim

\proclaim
{Theorem 3} ($0^\#$ does not exist) If II has a winning strategy in 
$\prg ,$ then II has a winning strategy in the descending chain game $\Cal G.$
\endproclaim
The descending chain game $\Cal G$ is to be defined shortly. An improvement
of the Theorem 3 appears in \cite {Vel}.

\proclaim
{Theorem 4} Cons($\kappa$ supercompact) implies 
Cons(in any $\aleph _1$-distributive algebra II has a winning 
strategy for $\prg ).$
\endproclaim

We also show that the latter statement has a consistency strength of 
at least that of $o(\kappa )=\kappa ^{++}.$

\demo
{Proof of the Theorem 1} Let $\sigma$ be a winning strategy for 
II in $\prg$ associated with a Boolean algebra $B.$ 
\proclaim
{Definition} If $\dot r$ is a $B$-name for a real, $j<\omega ,$ 
$\langle p,t\rangle \in B\times B^{\leq \omega },$ we say that 
$\langle p,t\rangle$ decides $\dot r(j)$ if for some $i<lth(t)$ 
$t(i)=\boldkey [\dot r(j)=1\boldkey ].$ If this is the case then we set 
$\langle p,t\rangle (\dot r(j))=1$ if $\sigma (\langle p,t\rangle )(i)=1$ 
and $0$ otherwise.
\endproclaim
 
Fix a sequence $\langle r_\alpha :\alpha <\omega _1\rangle$ of distinct 
reals (elements of $2^\omega )$ which does not contain a perfect subset. 
Let $\dot \alpha$ be a $B$-name for a countable ordinal.
\proclaim
{Claim 1} For each $\langle p,t\rangle \in B\times B^{<\omega }$ there 
is $\langle p,s\rangle \in B\times B^{<\omega }$ such that $t\subset s$ 
and for each $j<\omega ,$ $s_0,s_1\in B^{<\omega }$ if
 $s\subset s_0,s_1,$ $\langle p,s_0\rangle ,$ $\langle p,s_1\rangle$
 both decide $r_{\dot \alpha }(j)$ then $\langle p,s_0\rangle 
(r_{\dot \alpha }(j))= \langle p,s_1\rangle (r_{\dot \alpha }(j)).$
\endproclaim
\demo
{Proof of the Claim} By contradiction. (!\# \$ ?) Assume there is 
$\langle p,t\rangle$ witnessing the failure of the statement. By induction 
on $lth(\eta ),$ $\eta \in 2^{<\omega }$ we build $s_\eta \in B^{<\omega }$ 
so that \roster
\item $s_{\langle \rangle }=t,$ $\eta _0\subset \eta _1$ implies 
$s_{\eta _0}\subset s_{\eta _1}$
\item for all $\eta \in 2^{<\omega }$ there is $j,$ $lth(\eta )<j<\omega$ 
such that both $\langle p,s_{\eta ^\smallfrown 0}\rangle ,$
 $\langle p,s_{\eta ^\smallfrown 1}\rangle$ decide all the 
$r_{\dot \alpha }(i),$ $i\leq j$ and they decide 
$r_{\dot \alpha }(j)$ differently.
\endroster
There is no problem in the induction. Once we are done, for $x\in 2^\omega$
 we set $r_x:\omega \to 2,$ $r_x(i)=\langle p,\bigcup _{\eta \subset x}
s_\eta \rangle (r_{\dot \alpha }(i)).$ Since $\sigma$ is a winning strategy,

$$q=p\land \underset {\eta \subset x,\sigma (\langle p,s_\eta \rangle )
(i)=1}\to \bigwedge s_\eta (i)-\underset {\eta \subset x,\sigma 
(\langle p,s_\eta \rangle )(i)=0}\to \bigvee s_\eta (i)\neq 0$$

($q$ witnesses the success of $\sigma$ against $p,\bigcup _{\eta \
subset x}s_\eta .)$ Now obviously $q\Vdash r_{\dot \alpha }=r_x.$ 
So all the $r_x,$ $x\in 2^\omega$ are on the sequence $\langle 
r_\alpha :\alpha <\omega _1\rangle .$ However, the construction of 
$s_\eta$'s gives that $\{ r_x:x\in 2^\omega \}$ is a perfect set of 
reals, contradicting our assumption on 
$\langle r_\alpha :\alpha <\omega _1\rangle .$
\enddemo

Fix $N\prec H_\theta$ countable, with $B,\sigma ,\langle r_\alpha :
\alpha <\omega _1\rangle \in N,$ and choose $p\in N\cap B.$ Let 
$\langle \dot \alpha _i:i<\omega \rangle$ be the list of all 
$B$-names for countable ordinals in $N.$ Let 
$\omega =\bigcup _{j<\omega }Z_j$ be a disjoint union of infinite sets. 
We build $\langle t_n: n\in \omega \rangle$ by induction on $n$ so that 
\roster
\item $t_n\in B^{<\omega }\cap N,$ $t_0=\langle \rangle ,$ $n_0<n_1$ 
implies $t_{n_0}\subset t_{n_1}$
\item if $n$ is the first element of some $Z_j$ then let $t_{n+1}
\in B^{<\omega }\cap N$ be some finite sequence witnessing the Claim 1 
for $\dot \alpha _j,$ $\langle p,t_n\rangle$ in $N$ such that 
$\langle p,t_{n+1}\rangle $ decides $r_{\dot \alpha _j}(0)$
\item if $n$ is $i$-th element of some $Z_j,$ $i>1,$ 
then $t_{n+1}=t_n^\smallfrown \boldkey [r_{\dot \alpha _j}(i-1)=1\boldkey ].$
\endroster
Let $t=\bigcup _{n<\omega }t_n.$ Then

$$q=p\land \underset {\sigma (\langle p,t\restriction i+1\rangle )=1}
\to \bigwedge t(i)- \underset {\sigma (\langle p,t\restriction i+1
\rangle )=0}\to \bigvee t(i)\neq 0$$

since $q$ witnesses the success of $\sigma$ against $p,t.$ We claim 
that $q\Vdash N[G]\cap \omega _1=N\cap \omega _1,$ finishing the proof 
of the Theorem 1. For let $\dot \alpha _j\in N$ be a name for a 
countable ordinal. Let $n$ be the first element of $Z_j.$ Define 
$r:\omega \to 2$ by $r(i)=k$ if $\forall s\in B^{<\omega }$ if $t_{n+1}
\subset s$ and $\langle p,s\rangle$ decides $r_{\dot \alpha _j}(i)$ 
then $\langle p,s\rangle (r_{\dot \alpha _j}(i))=k.$ The choice of 
$t_{n+1}$ makes sure that $r$ is welldefined. Now all the parameters 
used in the definition of $r$ are in $N,$ so $r$ itself is in $N.$ 
Moreover, $q\Vdash r_{\dot \alpha _j}=r$ and so $q\Vdash \dot \alpha _j<N
\cap \omega _1.$ We are done.

\enddemo

This solves the Problem 2. For example shooting a closed unbounded set 
through a stationary costationary subset of $\omega _1$ \cite {BHK} does 
not add countable sequences of ordinals, but it is not semiproper and so 
there is no winning strategy for II in $\prg .$ In the end of the paper 
we give a more sophisticated example of a proper 
$\aleph _0$-distributive algebra with no winning strategy 
for II in $\prg$ (in ZFC).

\proclaim
{Definition} Given a complete Boolean algebra $B,$ the descending 
chain game $\Cal G$ is defined as follows:

\medpagebreak
\settabs \+\indent&abb&ccc&ddd&eee&fff&ratattatat&ggg&hhh&ratatttatat&&\cr
\+&I&$b_0$&&$b_1$&&\dotfill&$b_i$&&\dotfill&$i<\omega$&\cr
\+&II&&$c_0$&&$c_1$&\dotfill&&$c_i$&\dotfill&\cr
\medpagebreak

where $b_i,$ $c_i$'s are elements of $B$ and $b_0\geq c_0\geq b_1\geq c_1
\geq \dots .$ II wins a run of $\Cal G$ if $\underset {i<\omega }\to
\bigwedge b_i\neq 0.$
\endproclaim

\cite {Jech} has a proof that any $B$ for which II has a winning strategy 
in $\Cal G$ regularly embeds into an algebra with $\sigma$-closed dense 
subset. However, the question whether $B$ itself contains necessarily 
such a subset remains open.

\proclaim
{Fact 2}(\cite {Vel}) For any $B,$ II has a winning strategy in 
$\Cal G$ if and only if II has a winning strategy in the following game 
$\Cal G_1:$

\medpagebreak
\settabs \+\indent&abb&ccccc&ddd&eee&fff&ratattatat&ggg&hhh&ratatttatat&brum&&\cr
\+&I&$p,A_0$&&$A_1$&&\dotfill&$A_i$&&\dotfill&&$i<\omega$&\cr
\+&II&&$a_0$&&$a_1$&\dotfill&&$a_i$&\dotfill&\cr
\medpagebreak

where $p,a_i$'s are elements of $B,$ $A_i\subset B$ are maximal antichains 
and $a_i\in A_i.$ II wins a run of $\Cal G_1$ if $p\land 
\underset {i<\omega }\to \bigwedge a_i\neq 0.$
\endproclaim

(To pacify the reader, we state explicitly here that this is the 
last game introduced in this paper.)
\demo
{Proof of the Theorem 2} We follow closely the proof of the Theorem 1. 
Fix $B,\sigma$ such that $B$ has $(2^\omega )^+$-c.c. and $\sigma$ is a 
winning strategy for II in $\prg$ associated with $B.$ Due to the Fact 2, 
it is enough to find a winning strategy for II in $\Cal G_1.$ II proceeds 
as follows: at the $m$-th move he will have $p\in B,A_i :i\leq m,$ 
representing the moves of I, $a_i:i<m,$ II's answers and $\langle r_a^i:
i<m,a\in A_i\rangle ,$ reals $r_i,$ $i<m$ and $s_i\in B^{<\omega },$ 
$i\leq m$ so that \roster
\item $i<j\leq m$ implies $s_i\subset s_j$
\item $\langle r_a^i:a\in A_i\rangle$ is a sequence of distinct reals 
which does not contain a perfect set, for all $i<m.$
\item For $i<m,$ if $\dot a_i$ is the $B$-name for the unique element 
of $A_i$ in the generic ultrafilter, then the real $r_i$ is such that 
for all $j<\omega ,$ $s\in B^{<\omega },$ if $s_{i+1}\subset s,$ 
$\langle p,s\rangle$ decides $r_{\dot a_i}(j)$ then $\langle p,s\rangle 
(r_{\dot a_i}(j))=r_i(j).$ We have $r_i=r_{a_i}^i.$
\item $\langle p,s_m\rangle$ decides all $r_{\dot a_i},$ $i<m$ up to $m.$
\endroster
It is not difficult to proceed in the construction. Choose a sequence 
$\langle r_a^m:a\in A_m\rangle$ of distinct reals which does not contain 
a perfect subset. Here we use the chain condition of $B.$ Similarly as in 
the claim 1, we find an extension $s_{m+1}$ of $s_m$ and a real $r_m$ such 
that (3) continues to hold at $m.$ We can obviously require that 
$\langle p,s_{m+1}\rangle$ decides all $r_{\dot a_i}$ up to $m+1,$ 
$i\leq m.$ II plays the unique $a_m$ such that $r_{a_m}^m=r_i.$ 
This is a winning strategy since in the end,

$$
p\land \underset {i<\omega }\to \bigwedge a_i\geq p\land 
\underset {m<\omega ,\sigma (\langle p,s_m\rangle )(i)=1}\to 
\bigwedge s_m(i)-  \underset {m<\omega ,\sigma (\langle p,s_m\rangle )
(i)=0}\to \bigvee s_m(i)\neq 0
$$
the element of $B$ on the right hand side of the inequality witnessing 
the success of $\sigma$ against $p,\bigcup _ms_m.$ (Notice it decides 
all $r_{\dot a_i}$ to be equal to $r_i$ and so the $\dot a_i$'s to be 
equal to $a_i$'s.)

\enddemo

\proclaim
{Example 1} Theorem 2 is optimal; Cons($\kappa$ measurable) implies 
Cons( there is $B$ of density $\aleph _2$ where II wins $\prg$ but does not 
win $\Cal G).$
\endproclaim

Let $\kappa$ be measurable, $G\subset Coll(\omega _1, <\kappa )$ generic.
We work in $V[G].$ There is a variation of Namba forcing in which II
wins $\prg$ (see proof of the Theorem 4). We will need the following:
there is $C,$ a complete algebra in which II wins $\prg$ and there is a $C$
-name $\dot f,$ $C\Vdash$``$\dot f:\omega \to \omega _2^{V[G]}=\kappa$
is increasing and cofinal". Now fix $S\subset \{ \alpha <\omega _2:
cof(\alpha )=\omega \}$ such that both $S,$ $\{ \alpha <\omega _2:
cof(\alpha )=\omega \} \setminus S$ are stationary in $\omega _2$ and let
$P=\{ a\in [S]^{\aleph _0}:a$ closed $\}$ ordered by endextension. $B=RO(P).$
II does not win $\Cal G$ here since $B$ is not proper. We describe a winning
strategy for II in $\prg .$ Fix $\lessdot ,$ a wellordering of $P$ and
$\sigma ,$ a winning strategy for II in $\prg$ associated with $C.$
At $m$-th step of our game we will have $p,b_0,\dots b_m\in B,$
I's questions, $r_0,\dots r_{m-1}\in 2,$ II's answers and $C$-names
$\dot a_0,\dots \dot a_{m-1}$ so that \roster
\item $\forall i<m$ $C\Vdash$``$\dot a_i$ is $\lessdot$ least extension of
$\dot a_{i-1}$ (or $\check p$ if $i=0)$ deciding $b_i$ (that is, either
under $b_i$ or incompatible with $b_i)$ with $sup(\dot a_i)>\dot f(i)"$
\item in the game $\prg$ associated with $C,$ if I asks $1,\boldkey [\dot a_i
<b_i\boldkey ],$ $i<m,$ then at $i$-th position $\sigma$ answers just $r_i.$
\endroster
It is trivial for II to proceed in the construction. We want to show that
in the end, II wins. So let $p,b_0,b_1\dots ,b_i,\dots ,$ $i<\omega ,$ be
members of $B.$ Find $N\prec H_\theta$ countable such that $\lessdot ,p,
B,C,\sigma ,\langle b_i,i<\omega \rangle \in N$ and $sup(N\cap \omega _2)\in S.$
Find a $N$-generic filter $H\subset N\cap C$ such that $q\in H,$ where $q$
witnesses the success of $\sigma$ against $1, \boldkey [\dot a_i<b_i\boldkey ]$
$i<\omega .$ Let $a=\bigcup _{i<\omega }\dot a_i/H \cup \{ 
sup(N\cap \omega _2)\} .$ It is left to the reader to show that $a\in P$ and
that $a$ decides all the $b_i$'s as required.

\demo
{Proof of the Theorem 3} As in the Theorem 2, we show that a winning 
strategy for II in $\prg$ gives a winning strategy in $\Cal G_1$ and 
by Fact 2 we finish.
Thus let $\sigma$ be a winning strategy for II in $\prg$ associated 
with a Boolean algebra $B.$ We describe how II wins $\Cal G_1.$
 After $m$-th move of I in that game we will have $p\in B,A_i :
i\leq m,$ representing the moves of I, $a_i:i<m,$ the answers II has given 
so far and an auxiliary sequence $s_m\in B^{<\omega }$ such that 
\roster
\item $s_0=\langle \rangle ,$ $m_0<m_1\leq m$ implies $s_{m_0}\subset s_{m_1}$
\item for $i<m$ there is $j<lth(s_{i+1})$ such that $a_i$ is $s_{i+1}(j)$
 and $\sigma (\langle p,s_{i+1}\rangle )(j)=1.$
\endroster
We show how to obtain $s_{m+1},$ $a_m\in A_m$ so that (1),(2) 
continue to hold. Then II's $m$-th move in $\Cal G_1$ will be $a_m.$ If 
we succeed in doing that for all $m<\omega ,$ II wins the run of
 $\Cal G_1$ since 

$$
p\land \underset {i<\omega }\to \bigwedge a_i\geq p\land \underset 
{m<\omega ,\sigma (\langle p,s_m\rangle )(i)=1}\to \bigwedge s_m(i)-  
\underset {m<\omega ,\sigma (\langle p,s_m\rangle )(i)=0}\to \bigvee 
s_m(i)\neq 0
$$

since the element of $B$ standing on the right witnesses the success of 
$\sigma$ against $p,\bigcup _{m<\omega }s_m .$ 

To get the $m$-th move, II simulates a run of $\prg$ as follows: let $\theta$ 
be large regular, $\lessdot$ a wellordering of $H_\theta .$ Choose 
$N\prec \langle H_\theta ,\in ,\lessdot \rangle$ such that $s_m,A_m, p,
B\in N,$ $|N|=\aleph _0$ and there is no countable submodel $M$ of 
$\langle H_\theta ,\in ,\lessdot \rangle$ with $N\subset M,$ $B\cap 
(M\setminus N)\neq 0$ and $M\cap \omega _1= N\cap \omega _1.$ The 
existence of such a model is guaranteed by $\lnot 0^\#.$ Let 
$\{ b_i:i<\omega \} =N\cap A_m.$ Now we proceed as in the proof of 
the Theorem 1 to get a $N$-semimaster condition $q<p$ with two 
changes: $t_0=s_m$ and we require that $b_i$ sits on the sequence $t_{i+1}.$ 
We claim that $\exists i<\omega$ $\exists j<lth(t_{i+1})$ $b_i=t_{i+1}(j)$ 
and $\sigma (\langle p,t_{i+1}\rangle )(j)=1.$ ({\it Proof.} Otherwise 
$q<p,$ the condition witnessing success of $\sigma$ against 
$p,\bigcup _{i<\omega }t_i$ forces ``$N[G]\cap \omega _1=N\cap 
\omega _1,$ $A_m\cap N\cap G=0".$ Find $a\in A_m$ compatible with $q.$ 
$a\notin N$and a standard argument gives that if $M$ is the Skolem hull 
of $N\cup \{ a \}$ in $\langle H_\theta ,\in ,\lessdot \rangle$ using 
the canonical Skolem function generated by $\lessdot ,$ then $M\cap 
\omega _1=N\cap \omega _1.$ This would contradict maximality of $N.)$ 
So fix $i$ as above, set $s_{m+1}=t_{i+1}$ and let II play $a_m=b_i.$ 
The induction hypotheses continue to hold and we are done.
\enddemo

Obviously, the method of proof of the Theorem 3 falls short of 
reaching measurable cardinals. If we try to find a counterexample 
to its statement in Dodd's $K,$ the most natural thing to do is to 
look at projections of P in K, where $U$ is a normal measure and $P\in L[U]$ 
is the Prikry forcing. This will not work, though.

\proclaim
{Claim 2} Let $U$ be a normal measure at $\kappa ,$ $P$ the associated 
Prikry forcing. There are no nontrivial regular subalgebras $Q$ of $RO(P)$ 
with $Q\Vdash$``$\kappa$ is regular".
\endproclaim

\proclaim
{Corollary} In $L[U]:$ if $j:V\to M$ is the ultrapower by $U,$ then there 
are no regular subalgebras of the Prikry forcing in $M$ and so in $K.$
\endproclaim

{\it Remark.} In the presence of large cardinals it can happen that there is 
a regular subalgebra of a Prikry forcing in the associated ultrapower. Also,
 see Example 4 for a poset in $K$ with a winning strategy for II in $\prg$ 
in $L[U].$

\proclaim
{Question 1} Is it consistent to have $\kappa$ regular and a 
cardinal-preserving $Q$ of size $\kappa$ such that $Q$ singularizes $\kappa ?$
\endproclaim

In the view of the claim a negative answer to this question would prove that 
the Prikry forcing has no regular subalgebras of density $\kappa .$

\proclaim
{Question 2} Can the assumption of the Theorem 3 be weakened to ``$L[U]$ 
does not exist"? (Yes. \cite {Vel})
\endproclaim

\demo
{Proof of the Corollary} Work in $L[U].$ Let us assume that $Q\in M$ is 
an isomorphic copy of a dense subset of a regular subalgebra of 
the Prikry forcing. Then in $V,$ $Q\Vdash$``$\kappa$ is regular". If this
 were not the case then due to the $\kappa ^+$-c.c. of $Q$ in $V$ (and so 
in $M)$ the name witnessing that $Q\Vdash cof(\kappa )=\omega$ would be in 
$M.$ So $M\models$``$Q$ has $\kappa ^+$-c.c., does not collapse cardinals 
less than $\kappa$ and it singularizes $\kappa .$ So we have a failure 
of covering lemma in $M^Q$ \cite {DJ}. Therefore $Q\Vdash$``$\kappa$ is 
regular" and by the Claim 2 $Q$ is trivial.
\enddemo

\demo
{Proof of the Claim 2} Let us assume that we have $Q$ as in the Claim. 
First we prove that $Q$ does not add subsets of $\kappa$ and then argue 
that there are no $\kappa$-centered nontrivial forcings which do not add 
subsets of $\kappa .$ Since $Q,$ being a regular subalgebra of 
$\kappa$-centered $RO(P),$ is itself $\kappa$-centered, 
this will finish the proof.

Thus assume that $Q\Vdash$``$\dot A\in \Cal P(\kappa )\setminus V".$ 
$\dot A$ is a $P$-name and we build $A_t,t\in \kappa ^{<\omega }$ so 
that $\alpha \in A_t$ iff $\alpha \geq max(t)$ and there is $C\in U,$ 
$\langle t,C\rangle \Vdash _P\alpha \in \dot A.$ A standard argument 
yields $D\in U$ with $\langle 0,D\rangle \Vdash _P$``if $\kappa _0,
\kappa _1,\dots \kappa _i,\dots$ is the Prikry sequence then $\dot A\cap 
[\kappa _i,\kappa _{i+1})=  A_{\langle \kappa _0,\dots \kappa _i\rangle }
\cap [\kappa _i,\kappa _{i+1})."$ Choose $G\subset Q$ generic such 
that $proj_Q\langle 0,D\rangle \in G.$ Set $A=\dot A/G.$ Build 
$\lambda _0,\lambda _1,\dots \lambda _i ,\dots ,i<\omega$ in $V[G]$ by 
induction so that $\lambda _{i+1}=sup\{ \alpha <\kappa :\exists t\in 
\lambda _i^i$ $\alpha =min\{ \beta <\kappa :A\cap [max (t),\beta )\neq 
A_t\cap [max(t),\beta )\} \} .$ (Notice that since $A\notin V,$ for any 
$t\in \kappa ^{<\omega }$ there is $\alpha$ as in the definition of 
$\lambda _{i+1}.)$ Now it is easy to see that if $H\subset P,$ 
$\langle 0,D\rangle \in H,$ $G\subset H,$ is generic then the Prikry 
sequence has to be dominated by $\lambda _i$'s. However, $\kappa$ was 
regular in $V[G]$ and so $sup_i\lambda _i<\kappa ,$ a contradiction 
proving the Claim.
\enddemo

\demo
{Proof of the Theorem 4} Let $\kappa$ be supercompact, 
$G\subset Coll(\omega _1 ,<\kappa )$ generic. We claim that in 
$V[G]$ all $\aleph _1$-distributive algebras have winning strategies 
for II in $\prg .$ Work in $V[G].$ Let us assume that $B$ is 
$\aleph _1$-distributive, $\lambda =|\Cal P(B)|.$ In $V,$ pick a normal 
measure $U\subset \Cal P\Cal P_{<\kappa }\lambda .$ $U$ generates a filter 
$\hat U\subset \Cal P\Cal P_{\aleph _1 }\lambda$ in $V[G].$ It is possible 
to prove that $\Cal P\Cal P_{\aleph _1 }\lambda /I,$ where $I$ is the ideal
 dual to $\hat U,$ has a $\sigma$-closed dense subset $D$ (see \cite {Lav}).

Let $P=\{ T\subset (\Cal P_{\aleph _1 }\lambda )^{<\omega }:T$ is a tree, 
$T$ has a trunk $t$ and for $s\in T,$ $t\subset s$ implies $\{ a\in 
\Cal P_{\aleph _1 }\lambda :s^\smallfrown \langle a\rangle \in T\} 
\notin I\} ,$ ordered by inclusion. $P$ is a certain variation of Namba
 forcing. Now, $B$ regularly embeds into $RO(P):$ assume $H\subset P$ 
is generic. $H$ yields a sequence $\langle a_i:i<\omega \rangle \subset 
(\Cal P_{\aleph _1 }\lambda )^{V[G]}$ such that $\lambda =
\bigcup _{i<\omega }a_i.$ If $\langle D_\alpha :\alpha <\lambda \rangle$ 
is some enumeration in $V[G]$ of dense open subsets of $B,$ one can find 
a sequence $p_0>p_1>\dots >p_i>\dots ,i<\omega $ of elements of $B$ such 
that $p_i\in \bigcap _{\alpha \in a_i}D_\alpha ,$ which is open dense 
by distributivity of $B.$ Consequently, the filter on $B$ generated by 
the $p_i$'s will be generic over $V[G].$

To complete the proof, it is enough to show that II has a winning strategy 
in $\prg$ associated with $RO(P),$ since this property is obviously inherited 
by all of its subalgebras. We describe the strategy. At $m$-th step of the 
game we will have $T\in P,$ $b_i,i\leq m$ in $RO(P),$ the moves I made and 
an auxiliary sequence $T_i,i\leq m$ such that \roster
\item $T=T_0>T_1>\dots >T_m$ are trees with the same trunk $T,$
\item $T_{i+1}$ decides the statement ``$b_i\in H",$ $H$ generic, $i<m$
\item $\forall i\leq m$ $\forall s\in T_i$ $t\subset s$ implies $\{ a\in 
\Cal P_{\aleph _1 }\lambda :s^\smallfrown \langle a\rangle \in T_i\} \in D.$
\endroster
The construction of $T_{m+1}$ follows closely the proof of Prikry property 
in Prikry forcing. II answers 1 if $T_{m+1}\Vdash b_m\in H,$ otherwise 
the answer is 0. In the end, let $T_\omega =\bigcap _{i<\omega }T_i.$ By 
tree induction on elements of $T_\omega$ one can check that $T_\omega \in P$ 
is a tree with trunk $t.$ $\sigma$-closure of $D$ is used here. 
$T_\omega$ witnesses that II wins the run of $\prg$ 
in question. The Theorem follows.
 
\enddemo

\proclaim 
{Example 2} ($\lnot \exists \kappa$ $o(\kappa )=\kappa ^{++})$ A 
$\aleph _1$-distributive algebra for which II does not have a winning 
strategy in $\prg .$
\endproclaim

Let $\kappa$ be the $\omega _1$-st strong limit cardinal. Our hypothesis 
gives (see \cite {Git,Sch}) that $2^\kappa =\kappa ^+,$ $\kappa ^\omega =
\kappa ,$ $\square (\kappa ^+).$ Choose $E\subset cof(\omega )\cap 
\kappa ^+$ stationary, nonreflecting. The cardinal arithmetic implies 
$\lozenge (E).$ Fix $\langle S_\alpha :\alpha \in E\rangle$ so that 
$S_\alpha :\alpha ^{<\omega }\to 2^{<\omega }$ and for each 
$S:(\kappa ^+)^{<\omega }\to 2^{<\omega }$ there is $\alpha \in E$ 
such that $S\restriction \alpha ^{<\omega }=S_\alpha .$ Fix 
$\langle c_\alpha :\alpha \in E\rangle$ such that $c_\alpha \subset \alpha$ 
is cofinal of ordertype $\omega$ with an increasing enumeration 
$\gamma _0^\alpha ,\gamma _1^\alpha \dots .$ Define 
$\langle b_\alpha :\alpha \in E\rangle$ by $b_\alpha =\{ \gamma _i^\alpha 
\in c_\alpha :S_\alpha (\langle \gamma _0^\alpha ,\gamma _1^\alpha ,
\dots ,\gamma _i^\alpha \rangle )(i)=1\} .$ Finally we define our 
partial ordering $P=\{ a\subset \kappa ^+ :|a|\leq \kappa ,\forall 
\alpha \in E\cap sup(\alpha )+1$ $|(a\cap c_\alpha )\Delta b_\alpha |=
\aleph _0\} .$ Order by endextension.
\proclaim
{Claim 3} $\forall a\in P$ $\forall \alpha <\kappa ^+$ $\exists b<a$ 
$sup(b)\geq \alpha .$
\endproclaim
\demo
{Proof} Assume for contradiction that $a,\alpha$ witness the failure 
of the claim with $\alpha <\kappa ^+$ smallest possible.\roster
\item $\alpha$ is a successor. The contradiction is immediate. 
\item $\alpha \in E.$ $cof(\alpha )=\omega$ and one can build $a=a_0>a_1>
\dots >a_i >\dots ,i<\omega$ so that $sup(\bigcup _{i<\omega }a_i)=\alpha .$ 
Now we change $\bigcup _{i<\omega }a_i$ on $c_\alpha \setminus sup(a)$ to 
get $b$ so as to satisfy $|b \Delta b_\alpha |=\aleph _0.$ Since the 
ordertype of $c_\alpha$ is $\omega ,$ the similar statement at 
$\beta <\alpha$ will persist through this change. The amended set will 
therefore be a condition in $P$ endextending $a$ with supremum $\alpha .$
\item $\alpha \notin E,$ $cof(\alpha )=\omega .$ Build  $a=a_0>a_1>\dots >
a_i >\dots ,i<\omega$ so that $sup(\bigcup _{i<\omega }a_i)=\alpha .$ 
$\bigcup _{i<\omega }a_i\in P,$ it endextends $a$ and 
its supremum is $\alpha .$
\item $cof(\alpha )>\omega .$ Choose $C\subset \alpha ,$ $C\cap E=0.$ 
Build $a=a_0>a_1>\dots >a_i >\dots $ so that $j$ limit implies 
$sup(\bigcup _{i<j}a_i)\in C.$ Then $\bigcup _i a_i\in P,$ it endextends 
$a$ and its supremum is $\alpha .$
\endroster
In all cases, a contradiction. The Claim is proven.
\enddemo
So $P\Vdash$``$\bigcup G\subset \kappa ^+$ is cofinal". Let me prove that 
$P$ is $\aleph _1$-distributive. Let $\langle D_i:i<\omega _1\rangle$ be 
a sequence of open dense subsets of $P,$ $a\in P.$ Choose 
$M\prec H_\theta ,$ $\theta$ large, $\kappa ,P,a,$ $\langle D_i:
i<\omega _1\rangle \in M,$ $|M|=\kappa ,$ $M^\omega \subset M,$ $M\cap 
\kappa ^+=\alpha$ of cofinality $\omega _1 .$ Choose $C\subset \alpha$ 
closed unbounded of ordertype $\omega _1$ disjoint from $E.$ Much like 
in the case (4) in the proof of the Claim, we can construct a sequence 
of conditions $a>a_0>a_1>\dots a_i >\dots ,i<\omega _1$ such that each 
$a_i\in D_i\cap M,$ for $i<\omega _1$ limit $sup(\bigcup _{j<i}a_j)\in 
C.$ Then $\bigcup _{i<\omega _1}a_i\in \bigcap _{i<\omega _1}D_i$ and 
is less than $a,$ finishing the proof of distributivity. Now let $\sigma$ 
be a strategy for II in $\prg$ associated with $RO(P).$ Consider $S:
(\kappa ^+)^{<\omega }\to 2^{<\omega },$ the restriction of $\sigma$ 
to questions like $p=1,$ $\boldkey [\alpha \in \bigcup G \boldkey ],$ 
$\alpha <\kappa ^+.$ Fix $\alpha \in E$ such that $S_\alpha =
S\restriction \alpha ^{<\omega }.$ Then if I puts $p=1$ and $b_i=
\boldkey [\gamma _i^\alpha \in \bigcup G\boldkey ],$ then $\sigma$ must 
fail, because due to the construction of $P,$ no $a\in P$ with $sup(a)>
\alpha$ (and there is a dense set of those by the Claim 3) will agree 
with the choices $\sigma$ has made in the run of $\prg .$

\proclaim
{Example 3} (ZFC) A proper $\aleph _0$-distributive algebra for which 
II does not have a winning strategy in $\prg .$
\endproclaim

Let $\Cal T=\{ T: T$ is a tree on countable ordinals such that 
if $height(T)=\beta$ then $\forall x\in T$ $\forall \gamma <\beta$ 
$\exists y$ $<_T$ compatible with $x$ at the $\gamma$-th level of 
$T\} .$ For $T\in \Cal T,$ $\beta <\omega _1$ we define $T\restriction 
\beta =\{ x\in T:\exists \gamma <\beta$ $x$ is at level $\gamma$ of 
$T\} .$ We view $T\restriction \beta$ as a tree  with its ordering 
inherited from $T.$ Choose $f:\Cal T\to \Cal T$ such that $height(f(T))
= height(T)+1$ and $T=f(T)\restriction height(T).$ Finally we define a 
partial order $P=Q_0*\dot Q_1.$ Choose $S\subset \omega _1$ stationary
 costationary and set $Q_0=\{ T\in \Cal T:\forall \beta \in S\cap 
height(T)$ $T\restriction \beta +1=f(T\restriction \beta )\}$ with 
ordering $T_0<T_1$ if $T_1=T_0\restriction height(T_1).$ $Q_0$ 
is $\sigma$-closed. Let $Q_0\Vdash$``$\dot Q_1=\bigcup G$ with the 
obvious ordering". A standard argument gives that $Q_0\Vdash$``$\bigcup G$ 
is a Souslin tree" and so $P$ is $\aleph _0$-distributive proper, since it
 is an iteration of two such forcings. $P$ has a dense subset of size 
$2^\omega ,$ namely $\{ \langle T,x\rangle :T\in Q_0,x\in T\} .$ In 
view of the Theorem 2, to show that II has no winning strategy for 
$\prg$ associated with $RO(P),$ it is enough to show that there is 
no such strategy for $\Cal G.$ Let $   \sigma$ be a strategy for II 
in $\Cal G.$ W.l.o.g. the moves for II are from the abovedescribed 
dense subset of $P.$ Choose $N\prec H_\theta$ countable with $f,B,P,
\sigma \in N,$ $\beta =N\cap \omega _1\in S$ and $D_i: i<\omega$ an 
enumeration of open dense subsets of $Q_0$ in $N.$ Let $\langle 
\eta _i:i<\omega \rangle$ be an enumeration of $2^{<\omega }$ respecting 
the extension ordering of $\eta _i$'s. In $N$ we construct partial 
plays $W_i:i<\omega$ of $\Cal G$ by induction so that \roster
\item $\eta _i$ extends $\eta _j$ implies $W_i$ extends $W_j,$ 
$W_i$'s are played according to $\sigma$
\item $lth(W_i)=lth(\eta _i)$ and the last move in $W_i$ is played by II
\item $i<j$ implies that if $\langle T_i,x_i\rangle$ is the last move
 played in $W_i$ and $\langle T_j,x_j\rangle$ is the 
last move played in $W_j$ then $T_j<T_i$
\item if $\eta _i$ and $\eta _j$ are extension-incompatible and 
$i<j$ then $x_i,x_j$ are incompatible in $T_j$
\item $T_i\in D_i.$
\endroster
There is no problem in the induction. Once we are done, set $T=\bigcup T_i.$ 
If $x_{\eta _i}\in T$ is in the last move played by II in $W_i$ then
 the construction above makes sure that $\langle x_\eta :\eta \in 
2^{<\omega }\rangle$ form a perfect subtree of $T.$ So one of its 
branches does not have a lower bound in $f(T),$ which has only 
countably elements on $\beta$-th level. The game played along this 
branch witnesses the failure of $\sigma .$ (Notice that $f(T)$ is a 
``mandatory" extension of $T$ since $height(T)=\beta$ and $\beta \in S.)$

\proclaim
{Example 4} ($\kappa$ measurable) A proper algebra $B$ of density $\kappa$ 
where II wins $\prg$ but II does not win $\Cal G.$
\endproclaim
Let $S\subset \{ \alpha <\kappa : cof (\alpha )=\omega \}$ be such that  
both $S,$ $\{ \alpha <\kappa : cof (\alpha )=\omega \} \setminus S$ 
are stationary in $\kappa .$ Let $\Cal T,$ $f$ be as in Example 3. We define 
the forcing $P=Q_0*\dot Q_1$ as follows: $Q_0=\{ \langle g, T\rangle :g:
\alpha \to \kappa$ for some $\alpha <\omega _1$ is increasing continuous 
and $T\in \Cal T$ has height $\alpha$ and for all $\beta <\alpha ,$ if 
$g(\beta )\in S$ then $T\restriction \beta +1=f(T\restriction \beta ).$ 
Order by coordinatewise extension. $Q_0$ is $\sigma$-closed 
and $Q_0\Vdash$``$R=\bigcup \{ T\in \Cal T:T$ appears in a condition in 
the generic ultrafilter $\}$ is a Souslin tree". ({\it Proof.} Let
 $p\in Q_0$ force ``$\dot A\subset R$ is a maximal antichain". 
Choose $N\prec H_\theta$ countable with $S,p,\dot A\in N$ and 
$sup(N\cap \kappa )\notin S.$ Build a $N$-generic sequence 
$p>\langle g_0,T_0\rangle >\dots >\langle g_i, T_i\rangle >\dots ,
i<\omega .$ Then $C=\{ x\in \bigcup _iT_i:\exists i<\omega$
 $\langle g_i, T_i\rangle \Vdash \check x\in \dot A\}$ is a maximal 
antichain in $\bigcup _iT_i .$ Set $g=\bigcup _ig_i \cup \{ N\cap 
\kappa \}$ and find an extension $T$ of $\bigcup _iT_i$ of height 
$N\cap \omega _1 +1$ such that $C$ is still a maximal antichain in
 $T.$ Then $p>\langle g,T\rangle \Vdash$``$\dot A=\check C,$ in particular 
it is countable".) Let $Q_0\Vdash$``$\dot Q_1$ is $R$ with the
 obvious ordering". Let $B=RO(P).$ $B$ is proper and II has no winning 
strategy in $\Cal G,$ as one can show with a variation of the argument 
in Example 3. $B$ has a dense subset $D=\{ \langle g,T,x\rangle :
\langle g,T\rangle \in Q_0,$ $dom(f)=\alpha +1$ for some $\alpha 
<\omega _1$ and $x\in T$ is on the $\alpha$-th level of $T\} .$
 (For $p\in D,$ $p=\langle g,T,x\rangle$ we call $g(max(dom(g)))$ 
maximum of $p.)$ So $B$ has density $\kappa .$ To show that II 
has a winning strategy in $\prg ,$ proceed as follows: choose 
$\lessdot,$ a wellordering of $D$ and $U,$ a measure on $\kappa .$ 
At the $m$-th stage of the game we will have $p,b_0\dots ,b_m ,$ 
I's moves, $r_0,\dots r_{m-1},$ II's answers and an auxiliary 
sequence $A_i:i<m$ such that
\roster
\item $A_{m-1}\subset \dots \subset A_1\subset A_0$ are all members of $U$
\item $\forall i<m$ $\forall \alpha _0< \alpha _1<\dots <\alpha _i\in A_i$  
if $p>p_0>p_1>\dots >p_i$ is such that for $j\leq i$ $p_j$ is 
$\lessdot$-least extension of $p_{j-1}$ (or of $p$ if $j=0)$ 
deciding $b_j$ with  maximum of $p_j$ $\geq \alpha _ j$ then 
$p_j\Vdash b_j\in G$ if $r_j=1$ and $p_j\Vdash b_j\notin G$ if $r_j=0.$
\endroster
It is easy to proceed in the construction. Define a 
partition $h:A_{m-1}^{m+1}\to 2$ (or $h:\kappa \to 2$ if $m=0)$ by 
$h(\alpha _0,\dots \alpha _m)=1$ just in case when if $p>p_0>p_1>
\dots >p_m$ is such that for $j\leq m$ $p_j$ is $\lessdot$-least 
extension of $p_{j-1}$ (or of $p$ if $j=0)$ deciding $b_j$ with 
maximum of $p_j$ $\geq \alpha _ j$ then $a_m\Vdash b_m\in G.$ 
($\alpha _j$ are in increasing order.) There is $A\subset A_{m-1}$ 
in $U$ homogeneous for $h.$ Let $A_m=A$ and $r_m$ the homogeneous value.

Let us have an arbitrary run of $\prg$ according to the above strategy.
 We want to see that II wins. Let $A=\bigcap _{i<\omega }A_i\in U.$ 
Choose $N\prec H_\theta$ such that $p,b_i,i<\omega ,A,S,\lessdot \in N,$ 
$N\cap \kappa \notin S,$ $cof(N\cap \kappa )=\omega .$ Then one can find
 a sequence $\alpha _i:i<\omega \subset A\cap N$ increasing and cofinal 
in $N\cap \kappa .$ We construct $p>p_0>p_1>\dots p_i \dots ,i<\omega$ so 
that $\forall i<\omega$ $p_i$ is the $\lessdot$-least extension of $p_{i-1}$
 (or of $p,$ if $i=0)$ deciding $b_i$ with maximum of $p_i$ $\geq 
\alpha _i.$ Then since $N\cap \kappa \notin S$ one can find a 
condition extending all $p_i$'s and therefore witnessing the success of II.

The complexity of the previous example was due to our effort to make $B$ 
proper. The method shows that for example shooting a club subset through 
$S$ as above with countable approximations gives a Boolean algebra where 
II wins $\prg .$

\proclaim
{Question 3} The game $\prg$ was defined to have length $\omega .$ Are versions
of it with different countable lengths equivalent in terms of existence
of winning strategies?
\endproclaim

In view of the Theorem 3 one needs at least $0^\#$ for answering 
this question negatively. Example 4 could be relevant here, as the method
used to generate a winning strategy in $\prg$ gives no hint how to
proceed in longer games. Without further specifications on $\kappa$ and $S$
it could happen, however, that there would be winning strategies even for
longer games. For that, go to the Gitik's models \cite {Zapl} 
for Rudin-Keisler increasing
sequence of $\kappa$-complete ultrafilters at $\kappa$ of length $\alpha <
\omega _1$ and use similar trick as in the Example 1, only with forcing
$Q_\kappa$ instead of Namba forcing.


\Refs
\widestnumber \key {Dodd}
\ref
 \key BHK
 \by J. Baumgartner, L. Harrington, E. Kleinberg
 \paper Adding a closed unbounded set
 \jour Jour. Symb. Logic
 \vol 41
 \yr 1976
 \pages 481--482
\endref
\ref
 \key DJ
 \by A. J. Dodd, R. B. Jensen
 \paper The covering lemma for $L[U]$
 \jour Ann. Math. Logic
 \vol 22
 \yr 1982
 \pages 127--135
\endref
\ref
 \key Dodd
 \by A. J. Dodd
 \book The core model
 \bookinfo London Math. Soc. Lecture Notes in Mathematics no. 61
 \yr 1982
 \publaddr Cambridge
\endref
\ref
 \key Git
 \by M. Gitik
 \paper The strength of the failure of SCH
 \jour Jour. Pure Appl. Logic
 \vol 51
 \yr 1991
 \pages 215--240
\endref
\ref
 \key Jech
 \by T. J. Jech
 \paper More game-theoretic properties of Boolean algebras
 \jour Ann. Pure Appl. Logic
 \vol 26
 \yr 1984
 \pages 11--29
\endref
\ref
 \key Kan
 \by A. Kanamori
 \book Higher infinite
 \toappear
\endref
\ref
 \key Kun
 \by K. Kunen
 \book Set theory, an introduction to independence proofs
 \publaddr Amsterdam
 \yr 1980
\endref
\ref
\key Lav
 \by R. Laver
\endref
\ref
 \key P\v r\' \i
 \by K. L.  P\v r\' \i kr\' y
 \paper Changing measurable into accessible cardinals
 \jour Dissertationes Mathematicae
 \vol 68
 \yr 1970
 \pages 5--52
\endref
\ref
 \key Sch
 \by E. Schimmerling
 \book Combinatorial principles in the core model
 \bookinfo Ph. D. thesis, UCLA
 \yr 1992
\endref
\ref
 \key She
 \by S. Shelah
 \book Proper forcing
 \bookinfo Lecture Notes in Mathematics no. 940
 \yr 1982
 \publaddr Berlin
\endref
\ref
 \key Vel
 \by B. Veli\v ckovi\' c
 \paper Playful Boolean algebras
 \jour Transactions of AMS
 \vol 296
 \yr 1986
 \pages 727--740
\endref
\ref
 \key Zapl
 \by J. Zapletal
 \paper Gitik's forcing: an exposition
\endref
\endRefs
\enddocument